\documentclass{amsart}

\usepackage{amsfonts}
\usepackage{amsmath}
\usepackage{amsthm}
\usepackage{amssymb}
\usepackage[english]{babel}
\usepackage{graphics}
\usepackage{latexsym}
\usepackage{longtable} 
\usepackage{mathrsfs} 
\usepackage{graphicx}
\usepackage{wrapfig} 
\usepackage[pdftex,colorlinks=true]{hyperref}
\usepackage{enumitem}
\usepackage{esint}

\newtheorem{theorem}{Theorem}
\newtheorem{corollary}[theorem]{Corollary}
\newtheorem{lemma}[theorem]{Lemma}
\newtheorem{proposition}[theorem]{Proposition}

\newtheorem{claim}[theorem]{Claim}
\newtheorem{example}[theorem]{Example}

\theoremstyle{definition}
\newtheorem{definition}[theorem]{Definition}

\newcommand{\mL}{\mathcal{L}}

\newcommand{\mF}{\mathcal{F}}

\newcommand{\mO}{\mathcal{O}}

\newcommand{\B}{\mathcal{B}}

\newcommand{\E}{\mathrm{E}}
\newcommand{\D}{\mathrm{D}}
\renewcommand{\H}{\mathrm{H}}

\newcommand{\F}{\mathrm{F}}

\newcommand{\R}{\mathbb{R}}

\newcommand{\N}{\mathbb{N}}

\newcommand{\mB}{\mathbb{B}}

\newcommand{\noi}{\noindent}
\newcommand{\ms}{\medskip}

\newcommand{\Ga}{\Gamma}
\newcommand{\de}{\delta}
\newcommand{\De}{\Delta}
\newcommand{\e}{\varepsilon}
\newcommand{\si}{\sigma}

\newcommand{\Om}{\Omega}

\newcommand{\weak }{\, -\!\!\!\!-\!\!\!\rightharpoonup}

\newcommand{\larrow}{\longrightarrow}
\newcommand{\ot}{\otimes}

\newcommand{\ri}{\rightarrow}

\newcommand{\p}{\partial}
\newcommand{\sub}{\subseteq}
\newcommand{\set}{\setminus}
\newcommand{\by}{\times}

\newcommand{\tr}{\mathrm{tr}}
\DeclareMathOperator{\sgn}{sgn} 

\newcommand{\dist}{\mathrm{dist}}

\newcommand{\bt}{\begin{theorem}}\newcommand{\et}{\end{theorem}}
\newcommand{\bd}{\begin{definition}}\newcommand{\ed}{\end{definition}}
\newcommand{\bl}{\begin{lemma}}\newcommand{\el}{\end{lemma}}
\newcommand{\beq}{\begin{equation}}\newcommand{\eeq}{\end{equation}}
\newcommand{\bc}{\begin{claim}}\newcommand{\ec}{\end{claim}}
\newcommand{\bex}{\begin{example}}\newcommand{\eex}{\end{example}}
\newcommand{\bcor}{\begin{corollary}}\newcommand{\ecor}{\end{corollary}}
\newcommand{\bp}{\begin{proof}}\newcommand{\ep}{\end{proof}}

\newcommand{\BPL}{\medskip \noindent \textbf{Proof of Lemma} }

\newcommand{\BPP}{\medskip \noindent \textbf{Proof of Proposition} }
\newcommand{\BPT}{\medskip \noindent \textbf{Proof of Theorem} }

\numberwithin{equation}{section}

\begin{document}

\title[Second Order absolute minimisers]{Existence, Uniqueness and Structure of Second Order absolute minimisers}

\author{Nikos Katzourakis and Roger Moser}

\address{Department of Mathematics and Statistics, University of Reading, Whiteknights, PO Box 220, Reading RG6 6AX, United Kingdom}
\email{n.katzourakis@reading.ac.uk}

\address{Department of Mathematical Sciences, University of Bath, Claverton Down, 
Bath BA2 7AY, United Kingdom}

  \thanks{\!\!\!\!\!\!\texttt{N.K. has been partially financially supported by the EPSRC grant EP/N017412/1}}
  
  \email{r.moser@bath.ac.uk}

\date{}

\keywords{$\infty$-Laplacian; $\infty$-Bilaplacian; Second Order absolute minimisers; Calculus of Variations in $L^\infty$.}

\begin{abstract} Let $\Omega \subseteq \mathbb{R}^n$ be a bounded open $C^{1,1}$ set. In this paper we prove the existence of a unique second order absolute minimiser $u_\infty$ of the functional
\[
\ \ \ \ \ \ \mathrm{E}_\infty (u,\mathcal{O})\, :=\, \| \mathrm{F}(\cdot, \Delta u) \|_{L^\infty( \mathcal{O} )}, \ \ \ \mathcal{O} \subseteq \Omega \text{ measurable},
\]
with prescribed boundary conditions for $u$ and $\mathrm{D} u$ on $\partial \Omega$
and under natural assumptions on $\F$. We also show that $u_\infty$ is partially smooth and there exists a harmonic function $f_\infty \in L^1(\Om)$ such that
\[
\F(x, \Delta u_\infty(x)) \, =\, e_\infty\, \mathrm{sgn}\big(f_\infty(x)\big)
\]
for all $x \in \{f_\infty \neq 0\}$, where $e_\infty$ is the infimum of the global energy.
\end{abstract}

\maketitle


\section{Introduction} \label{section1}

For $n\in \N$, let  $\Om \sub \R^n$ be a bounded open set and let also 
$\F : \Omega \times \R \larrow \R$ be a real function that is $\mL(\Om) \ot \B(\R)$-measurable, namely, measurable with respect to the product $\si$-algebra of the Lebesgue subsets of $\Om$ with the Borel subsets of $\R$. In this paper we consider variational problems for second order supremal functionals of the form
\beq \label{1.1}
\ \ \ \ \ \ \mathrm{E}_\infty (u,\mathcal{O})\, :=\,\big \| \mathrm{F}(\cdot, \Delta u) \big\|_{L^\infty( \mathcal{O} )}, \ \ \ \mathcal{O} \subseteq \Omega \text{ measurable},
\eeq
where the admissible functions $u$ range over the (Fr\'echet) Sobolev space
\beq \label{1.2}
\mathcal{W}^{2,\infty}(\Om)\,:= \bigcap_{1<p<\infty} \Big\{ u \in W^{2,p}(\Om)\ :\ \De u \in L^\infty(\Om)\Big\}.
\eeq
This is the natural space for our problem, and while it is clear that $W^{2, \infty}(\Omega) \subseteq \mathcal{W}^{2, \infty}(\Omega)$, the inclusion is strict in general if $n\geq2$. 

We defer a discussion of the wider context of \eqref{1.1}--\eqref{1.2}, but note that a central feature of supremal functionals is that global minimisers are typically not the appropriate ``minimal" objects because they usually fail to be optimal on subdomains. The latter is a property automatic for integral functionals, due to their additivity in the domain argument. Instead, we will consider the following localised variational notion.

\bd[Second order absolute minimisers]  \label{definition1}
A function $u\in \mathcal{W}^{2,\infty}(\Om)$ is called a second order absolute minimiser of \eqref{1.1} on $\Om$ if
\[
\E_\infty (u,\mO)\, \leq \, \E_\infty(u+\phi,\mO)
\]
for all open sets $\mO \sub \Om$ and all $\phi \in \mathcal{W}^{2,\infty}_0(\mO)$.
\ed
Here we have used the obvious notation $\mathcal{W}_0^{2,\infty}(\mO):= \mathcal{W}^{2,\infty}(\mO) \cap W_0^{2,2}(\mO)$. Given $u_0 \in \mathcal{W}^{2, \infty}(\Om)$, we will also write $\mathcal{W}_{u_0}^{2, \infty}(\Om) = u_0 + \mathcal{W}_0^{2, \infty}(\Om)$. The main goal in this paper is to establish the existence and uniqueness of second order absolute minimisers of \eqref{1.1} given fixed Dirichlet boundary conditions on $\p\Om$, deriving also additional information for the structure of such objects, including their optimal regularity. 

Although a global minimiser of $\E_\infty(\cdot ,\Om)$ may not always be a second order absolute minimiser (see e.g.\ Example \ref{ex5}), if $u \in \mathcal{W}^{2,\infty}(\Om)$ minimises $\E_\infty(\cdot, \Om)$ \emph{uniquely} in $\mathcal{W}^{2, \infty}(\Om)$ with respect to its own boundary conditions, then in fact $u$ is the unique second order absolute minimiser. Accordingly, we first give a condition guaranteeing that $u$ is the unique minimiser for its boundary values on any subdomain. While it is not clear a priori that this condition can be met, we subsequently prove that it is satisfied by exactly one function under given boundary conditions and very mild additional assumptions.

In the following, we will use the symbolisation ``$\sgn$" for the sign function, with the convention that $\sgn(0) = 0$. We will assume also that

\beq \label{as1}
\text{$\Om$ is a bounded connected open subset of $\R^n$, $n\geq 1$,} 
\eeq
and
\beq \label{as2}
\ \
\left\{\ \
\begin{split}
& \text{the function $\F : \Om \by \R \larrow \R$ is $\mL(\Om) \ot \B(\R)$-measurable and}
\\
& \text{for a.e. $x \in \Om$, $\xi \mapsto \F(x,\xi)$ is strictly increasing with $\F(x, 0) = 0$.} 
\end{split}
\right.
\eeq
Our first main result therefore is:

\bt[Criterion for unique minimisers] \label{uniqueness}
Suppose \eqref{as1}--\eqref{as2} hold and consider \eqref{1.1} and a function $u_* \in \mathcal{W}^{2,\infty}(\Om)$. If there exist a number $e_* \geq 0$ and a function $f_* \in L^1(\Om)$ satisfying
\beq \label{harmonic}
\Delta f_* \,=\, 0, \ \ \ \text{on $\Om$},
\eeq
such that
\beq \label{representation}
\ \ \F(\cdot, \Delta u_*) \,=\, e_*\sgn(f_*), \quad \text{a.e.\ on $\Omega$},
\eeq
then
\[
\E_\infty(u_*,\mO)\, <\, \E_\infty\big(u_* + \phi, \mO\big) \smallskip
\]
for any open $\mO\sub \Om$ and any $\phi \in \mathcal{W}_0^{2, \infty}(\mO) \set\{0\}$. 
\smallskip

Namely, for any open subset $\mO\sub \Om$, the function $u_*$ is the unique global minimiser of $\E_\infty(\cdot, \mO)$ in $\mathcal{W}_{u_*}^{2, \infty}(\mO)$. 
\et

Clearly, \eqref{representation} implies that $e_* = \E_\infty(u_*,\Om)$, unless $f_* \equiv 0$.
If we fix boundary data $u_0 \in \mathcal{W}^{2, \infty}(\Omega)$, then there exists at most one minimiser $u_* \in \mathcal{W}^{2,\infty}_{u_0}(\Om)$ and $e_*$ is uniquely determined as the infimum of $\E_\infty(\cdot,\Om)$ over the space. However, \emph{$f_*$ cannot be determined uniquely by the boundary conditions}. Note further that by assumption \eqref{as2}, the equation \eqref{representation} is equivalent to the \emph{representation formula}
\[
\ \ \ \ \ \ \Delta u_*(x) \,=\, \F(x,\cdot)^{-1} \Big(e_* \sgn\big(f_*(x)\big)\Big), \quad \text{ a.e. }x\in \Om.
\]
Moreover, by using standard argument involving Green functions (see e.g.\ \cite[Ch.\ 2]{GT}), we could represent $u_*$ in terms of $\F,f_*,e_*,u_* |_{\p\Om},\D_\nu u_*|_{\p\Om}$.

Theorem \ref{uniqueness} gives a connection between the variational problem and a PDE system of second order equations with a parameter consisting of \eqref{harmonic} and \eqref{representation}, which \emph{we may think of as a PDE formulation of the $L^\infty$ variational problem}. There exists, however, a more conventional analogue of the ``Euler-Lagrange equation" for \eqref{1.1}. This is the fully nonlinear PDE of third order
\beq \label{1.4}
\ \ \ \F(\cdot,\De u)\F_{\xi}(\cdot,\De u)  \big|\D \big(|\F(\cdot,\De u)|^2\big) \big|^2 \,=\, 0 \ \ \text{ on }\Om,
\eeq
where the subscript $\F_\xi$ denotes the partial derivative in $\xi$. The particular case arising from $\F(x,\xi)=\xi$ is coined the \emph{$\infty$-Bilaplacian}.  Since we will not make direct use of \eqref{1.4}, we refer to the recent paper \cite{KP2} wherein it is studied directly in a general context.

We will see now that under certain assumptions on $\Om$, $\F$ and the boundary data, the system \eqref{harmonic}--\eqref{representation} has in fact a solution $(u_*,f_*,e_*)$, with $f_* \not\equiv 0$ if $e_* > 0$. To this end, we need to assume additionally that $\F : \Om\by \R \larrow \R$ satisfies
\beq \label{1.3}
\ \ \left\{ \ \ 
\begin{split}
 & \ \ \ \ \ \ \ \F \in \,  C^2(  \Omega \times \R),
 \\
 & \ \  \  \ \ \F(x,0)\, =\, 0, \ \ \ x \in \Om,
\\
\exists \   c> 0   \ :  \ 
& \left\{
\begin{array}{ll}
c\,\leq\, \F_\xi(x,\xi) \, \leq\, \dfrac{1}{c} , &  (x,\xi)  \in \Om \by \R, \smallskip
\\
\F(x,\xi)  \, \F_{\xi\xi}(x,\xi)  \, \geq\, - \dfrac{1}{c} , &   (x,\xi)  \in  \Om \by \R.
\end{array}
\right.
\end{split}
\right.
\eeq
These conditions may seem restrictive at first, but since the solutions of the variational problem would be
the same for $g \circ F$ for any continuous, increasing, odd function $g \colon \R \to \R$, they do effectively allow
a wide range of functions. Condition \eqref{1.3} does imply that for any fixed $x\in \Om$, the function $|\F(x,\cdot)|$ is level-convex on $\R$ (i.e.\ has convex sublevel sets) but it does not require convexity. In addition to existence, we will prove that the absolute minimisers of the $L^\infty$ problem can be approximated by corresponding minimisers of $L^p$ problems as $p \to \infty$, obtaining also additional partial regularity information. For $p \in (1, \infty)$, we therefore define
\[
\ \ \ \ \E_p(u) \,:=\, \left( \fint_\Om \big| \F(\cdot,\De u)\big|^p\right)^{1/p}, \ \ \ u\in W^{2,p}(\Om),
\]
where the slashed integral sign denotes the average over $\Om$.

\bt[Existence, structure and approximation] 
\label{existence_theorem}
Let $\Om$ satisfy \eqref{as1} and also have $C^{1,1}$ boundary. Let also $\F$ satisfy \eqref{1.3} and fix a function $u_0$ in $\mathcal{W}^{2,\infty}(\Om)$ with $\De u_0 \in C(\overline{\Om})$.

\begin{enumerate}[label={\rm (\Roman*)}] 

\item \label{(I)} There exist a global minimiser $u_\infty \in \mathcal{W}_{u_0}^{2,\infty}(\Om)$ of $\E_\infty(\cdot,\Om)$ and a harmonic function $f_\infty \in L^1(\Omega)$ such that
\[
\ \ \ \F(\cdot, \Delta u_\infty) \, =\, e_\infty \sgn(f_\infty), \quad \text{a.e.\ on $\Om$},
\]
where $e_\infty = \E_\infty(u_\infty, \Om)$. Further, $f_\infty \not\equiv 0$ if $e_\infty >0$.

\item \label{(II)} Let
\[
\Gamma_\infty \, := \, f_\infty^{-1}(\{0\}). 
\]
If $e_\infty>0$, then $u_\infty$ belongs to $C^{3,\alpha}(\Omega \setminus \Gamma_\infty)$ for any $\alpha \in (0, 1)$ and $\Ga_\infty$ is a Lebesgue nullset. If $e_\infty=0$, then $u_\infty$ is a harmonic function.

\item \label{(III)} For $p \in \N$, let 
\[
\ \ \ e_p\, :=\, \inf\big\{ \E_p(u)\, : \, u \in W^{2,p}_{u_0}(\Om)\big\}.
\]
Then, for any $p$ large enough there exists a global minimiser $u_p \in W_{u_0}^{2, p}(\Om)$ of $\E_p$ satisfying $\E_p(u_p) = e_p$. Moreover, $e_{p} \larrow e_\infty$ as $p\ri\infty$ and there exists a subsequence $(p_\ell)_{\ell =1}^\infty$ such that $u_{p_\ell} \weak u_\infty$ in the weak topology of $\mathcal{W}^{2,\infty}(\Om)$ as $\ell \ri\infty$. In addition,
\[
\left\{
\begin{array}{ll}
u_{p_\ell} \larrow u_\infty, & \text{ in }C^{1}(\overline{\Om}), \smallskip
\\
\D^2 u_{p_\ell} \weak \D^2 u_\infty, & \text{ in }L^q(\Om,\R^{n\by n}) \text{ for all }q \in (1,\infty), \smallskip
\\
\De u_{p_\ell} \larrow \De u_\infty, &  \text{ a.e.\ on $\Om$ and in $L^q(\Om)$ for all $q \in (1,\infty)$},
\end{array}
\right.
\]
as $\ell \ri\infty$. Furthermore, $\De u_{p_\ell} \larrow \De u_\infty$ locally uniformly on $\Omega \setminus \Gamma_\infty$ if $e_\infty >0$ and locally uniformly on $\Omega$ if $e_\infty =0$.

\item \label{(IV)} We set
\beq \label{def_f_p}
f_p\, :=\, \frac{1}{e_p^{p - 1}}\big| \F(\cdot,\De u_p)\big|^{p-2}\, \F(\cdot,\De u_p)\,\F_\xi (\cdot,\De u_p)
\eeq
if $e_p \not= 0$, and $f_p := 0$ if $e_p = 0$.
Then, the harmonic function $f_\infty$ in \ref{(I)} may be chosen such that $f_{p_\ell} \larrow f_\infty$ as $\ell \ri \infty$ in the strong local topology of $C^\infty(\Omega)$.

\end{enumerate}
\et

By invoking Theorem \ref{uniqueness}, an immediate consequence is that the modes of convergence in Theorem \ref{existence_theorem}\ref{(III)} as $p\ri \infty$ are actually \emph{full and not just subsequential}. Also, known results on the regularity of nodal sets of solutions to elliptic equations \cite{HS} imply that \emph{$\Gamma_\infty$ is countably rectifiable}, being equal to the union of countably many smooth $(n - 1)$-dimensional submanifolds of $\Om$ and a set of dimension $(n-2)$ or less.

In most cases, the singular set $\Gamma_\infty$ is non-empty and divides $\Omega$ into two distinct parts, forcing the Laplacian of $u_\infty$ to have a jump across $\Gamma_\infty$ and no higher global regularity than membership in $\mathcal{W}^{2, \infty}(\Om)$ can be expected. The numerical and explicit solutions in \cite{KP2} confirm this, even when $n=1$. Further, \ref{(I)}--\ref{(IV)} above have been obtained therein for $n=1$ and in some other special cases (although have not been stated in this explicit fashion), whilst the qualitative behaviour emerging here was also observed numerically on the plane. The condition $\Delta u_0 \in C(\overline{\Omega})$ on the boundary data is required for our proof. The problem of existence is completely open without this assumption. 

When combined, Theorems \ref{uniqueness} and \ref{existence_theorem} imply in particular the following. 

\bcor \label{corollary4} Under the hypotheses of Theorem \ref{existence_theorem}, there exists a unique global minimiser $u_\infty$ of $\E_\infty(\cdot,\Om)$ in $\mathcal{W}_{u_0}^{2, \infty}(\Om)$, which is a second order absolute minimiser and a strong solution to the Dirichlet problem for \eqref{1.4}:
\[
\left\{
\begin{array}{rl}
\F(\cdot,\De u)\F_{\xi}(\cdot,\De u)  \big|\D \big(|\F(\cdot,\De u)|^2\big) \big|^2 \,=\, 0,  &\ \ \text{ in }\Om,
\\
u \,=\, u_0,  &\ \ \text{ on }\p \Om,
\\ 
\D u \,=\, \D u_0, &\ \ \text{ on }\p \Om,
\end{array}
\right.
\]
More precisely, $u_\infty$ is thrice differentiable a.e.\ on $\Om$ and satisfies the PDE in the pointwise sense.
\ecor

We recall that the fully nonlinear PDE above has essentially been derived from $L^p$-approximate equations and studied in the paper \cite{KP2} (see also \cite{A4, AB}). Note also that in contrast to the system \eqref{harmonic}, \eqref{representation}, equation \eqref{1.4} does not give uniqueness of strong solutions \emph{without additional supplementing conditions}. 
For example, for $\Omega = (-2, 2) \subseteq \R$ and $F(x, \xi) = \xi$,
consider
\[
u(x) = \begin{cases}
(x + 2)^2 - 3, & \text{if $-2 < x < -1$}, \\
2x, & \text{if $-1 \le x \le 1$}, \\
3 - (x - 2)^2, & \text{if $1 < x < 2$}.
\end{cases}
\]
This function satisfies \eqref{1.4} everywhere except at $\pm 1$ and is in $W^{2,\infty}(-2,2)$. But Theorems \ref{uniqueness} and \ref{existence_theorem} imply that there exists another strong solution $v$ with the same boundary data such that $v''$ takes one of only two values almost everywhere. We note that general explicit solutions in the case of $n=1$ (and no $x$-dependence) are given in \cite{KP2}.


We now connect the $L^\infty$ problem we study herein to the wider context. The study of supremal functionals and of their associated equations is known as Calculus of Variations in the space $L^\infty$. Variational problems for first order functionals of the form 
\beq \label{1.5}
\ \ (u,\mO) \mapsto\, \underset{x\in \mathcal{O}}{\mathrm{ess}\,\sup}\, \H \big(x, u(x),\mathrm{D} u(x)\big) ,\ \ \ u\in W^{1,\infty}(\Om),\ \mO \sub \Omega\text{ measurable},
\eeq
together with the associated PDEs, first arose in the work of Aronsson in the 1960s \cite{A1}--\cite{A3}. The first order case is very well developed and the relevant bibliography is very  extensive. For a pedagogical introduction to the theme which is accessible to non-experts, we refer to the monograph \cite{K8} (see also \cite{ACJ, C, C2}). This theory makes extensive use of the notion of viscosity solutions of
Crandall-Ishii-Lions (see \cite{C, CIL, K8}) for the underlying Aronsson
equation, which is a second order PDE. Without any attempt to be exhaustive, some milestones in the development of the theory evidencing the great interest the problem has attracted, can be found e.g.\ in \cite{ACJS, AS, BJW1, BJW2, BDM, CEG, CWY, ES, ESm, J, S2, WY, Y}. However, viscosity solutions fail for higher order equations, and therefore it is no surprise that we need to take a different approach for our problem. The approximation of the $L^\infty$-problem by $L^p$-problems, on the other hand, is very common in the first order theory, too.

Second order variational problems in $L^\infty$ have only very recently begun to be investigated and are still poorly understood. In the very recent paper \cite{KP2} (see also \cite{KP3}), the study of higher order variational problems and of their associated PDEs was initiated, focusing on functionals of the form
\beq \label{Polylap}
\ \ \ \ \   (u,\mO) \mapsto \, {\mathrm{ess}\, \sup}_\mathcal{O} \H \big(\mathrm{D}^2 u\big) ,\ \ \ u\in W^{2,\infty}(\Om),\ \mO \sub \Omega \text{ measurable}.
\eeq
Preliminary investigations on the second order $1$-dimensional case $\H(\cdot, u, u', u'')$ had previously been performed via different methods by Aronsson and Aronsson-Barron in \cite{A4, AB}. Very recently, arguments inspired by the present paper have been applied in \cite{KPa} to eigenvalue problems for the $\infty$-Bilaplacian.


It is remarkable that for our specific variational problem, \emph{global minimisers are
unique and are automatically absolute minimisers}. This is in contrast to first order problems in the Calculus of Variations in $L^\infty$ of the form
\[
\ \ \ \ \   (u,\mO) \mapsto \, {\mathrm{ess}\,\sup}_\mathcal{O} \H \big(\cdot,u,\mathrm{D} u\big) ,\ \ \ u\in W^{1,\infty}(\Om),\ \mO \sub \Omega \text{ measurable}
\]
which are very different. For these, global minimisers are neither unique nor absolute in general. In fact, the only case that a global minimiser of \eqref{1.1} can fail to be absolute is on a disconnected domain. In this sense, our Example \ref{ex5} given below is optimal. This unusual structure is a consequence of the properties of the Laplace operator and is not true for the full Hessian case of \eqref{Polylap}. In addition, for our problem, we obtain a fair amount of detailed information about the structures of the solutions with relatively simple means. Our methods are almost exclusively variational and essentially do not utilise any deep PDE machinery, in contrast to the seminal paper \cite{J}, wherein uniqueness of absolute minimisers was established for the first time --in the general case-- by using PDE methods. However, our techniques take advantage of the special structure of the functional and may be suitable only for specific problems. Nonetheless, some of the techniques that underpin Theorem \ref{existence_theorem} have been successfully deployed in \cite{MS, S1} to problems somewhat different to \eqref{1.1} (with $F$ having $u$-dependence), which suggests that further generalisation might be possible. In order to keep the presentation simple, however, we do not explore this direction any further in this work.

Apart from the intrinsic mathematical interest, the motivation to study higher order $L^\infty$ minimisation problems comes from several diverse areas. Minimisation problems in $L^\infty$ similar to the above (but with additional lower order terms) have also been studied in the context of differential geometry and in relation to the Yamabe problem in \cite{MS, S1}. Also, a prevalent applied problem in Data Assimilation (e.g.\ geosciences) and PDE-constrained optimisation (e.g.\ aeronautics) is the construction of approximate solutions to second order ill-posed PDE problems, see the model problem in \cite{K9} and references therein, as well as the classical monograph \cite{L}. For instance, in the modelling of aquifers, one needs to solve $\De u =f$ coupled with a pointwise constraint $K(u)=k$ for given $f, K, k $. By minimising the error $e(u):=|\De u -f |^2 + |K(u)-k|^2$ in $L^\infty$, one obtains uniformly (globally and even absolutely) best approximations. 

Note that the constrained problem ``$\De u =f$ \& $K(u)=k$" is generally overdetermined already, even without imposing any boundary conditions (either zeroth order or first order). In applications, one is typically interested in a variety of additional boundary conditions, including the case of partial or no boundary conditions. The methods developed in this paper do not directly extend to the case of (absolutely) minimising $\|e(u)\|_{L^\infty}$ and investigation of this type of problems is left for future work. Although global minimisers are easy to find, the interplay between global and absolute minimisers becomes very delicate due to the lower order terms. On the other hand, the variational method of approximate solutions allows significant flexibility, with the case of no boundary conditions corresponding to free minimisation in $\mathcal{W}^{2,\infty}$. The degrees of freedom in our homogeneous higher order problem require us to impose first order boundary conditions to obtain uniqueness.


\ms

\section{Proofs} \label{section2}

In this section we establish the proofs of Theorems \ref{uniqueness} and \ref{existence_theorem} and of Corollary \ref{corollary4}. Before delving into this, we give a simple example confirming that global minimisers of second order supremal functionals in general do not minimise on subdomains.

\begin{example}\label{ex5}{\rm For $\Om = (-1,0)\cup (0,1)$, consider the functional $\E_\infty(u,\mO)=\|u''\|_{L^\infty(\mO)}$ and let  $Q(x):=x \chi_{(-1,0)}(x) -x(x-1)\, \chi_{(0,1)}(x)$, $x \in \Om$. Then $Q$ is a global minimiser of $\E_\infty(\cdot,\Omega)$ in $\mathcal{W}_Q^{2,\infty}(\Om)$ with $\E_\infty(Q,\Om)=2$ and is also a second order absolute minimiser. However, for any function $\zeta \in C^\infty_c(-1,0)$ with $0<\|\zeta''\|_{L^\infty(-1,0)}<1$, the perturbation $Q+\zeta$ still satisfies $\E_\infty(Q+\zeta,\Om)=2$ and lies in $\mathcal{W}_Q^{2,\infty}(\Om)$, but does not minimise $\E_\infty(\cdot,(-1,0))$ over $W^{2,\infty}_Q((-1,0))$ because the only minimiser on $(-1,0)$ with boundary data $Q$ is the identity.}
\end{example}

We now continue with the proof of our first main result.

\BPT \ref{uniqueness}. Fix $\phi \in \mathcal{W}^{2,\infty}_{0}(\Om)$ with $\phi \not \equiv 0$ on $\Om$. Since $f_*$ is a harmonic function in $L^1(\Om)$, it follows that
\beq \label{2.18}
\int_\Om f_* \, \De \phi \, =\, 0.
\eeq
We set
\[
\Gamma_* \,:= \, f_*^{-1}(\{0\}).
\]
By standard results on the nodal set of solutions to elliptic equations \cite{HS} and the connectedness of $\Om$, it follows that if $f_* \not\equiv 0$ then $\Gamma_*$ is a Lebesgue nullset and if $f_* \equiv 0$ then $\Gamma_*=\Om$.

Let us first consider the case $f_* \not\equiv 0$. Note that $\De \phi$ cannot vanish almost everywhere on $\Om$ (as this would imply that $\phi \equiv 0$ by uniqueness of solutions of the Dirichlet problem for the Laplace equation). Therefore, we deduce that $f \De \phi \neq 0$ on a subset of positive Lebesgue measure in $\Om$. Hence, \eqref{2.18} implies that there exist measurable sets $\Om^\pm \sub \Om$ with $\mL^n(\Om^\pm)>0$ such that 
\[
\pm  f_*\,  \De \phi \, >\, 0, \ \ \ \text{a.e.\ on }\Om^\pm,
\]
where $\mL^n$ denotes the $n$-dimensional Lebesgue measure.  If \eqref{representation} holds true, then we have $|\F(\cdot, \De u_*)| = e_*$ a.e.\ on $\Om$ and
\[
\sgn(\De u_*) \,= \, \sgn(f_*) \, =\, \sgn(\De \phi), \ \ \ \text{ a.e.\ on }\Om^+. 
\]
As $\F(x, \cdot)$ is strictly increasing for a.e.\ $x \in \Om$, this gives that 
\[
|\F(\cdot, \De u_* + \De \phi)| \, >\, |\F(\cdot, \De u_*)| \,=\, e_*,
\ \ \text{ a.e.\ in }\Om^+. 
\]
Therefore,
\[
\E_\infty(u_*, \Om) \,=\, e_* \,<\, \E_\infty(u_* + \phi, \Om).
\]
It remains to consider the case $f_* \equiv 0$. Then, the hypothesis \eqref{representation} implies that $\Delta u_* = 0$ almost everywhere and so $\E_\infty(u_*,\Om) = 0$. On the other hand, by arguing as above it follows that $\E_\infty(u_* + \phi, \Om) >0$ for any $\phi \in \mathcal{W}_0^{2, \infty}(\Om) \setminus \{0\}$. Hence, we arrive at the same conclusion.

Finally, if $\mO \sub \Om$ is a non-empty open set and $\phi \in \mathcal{W}_0^{2, \infty}(\mO)\set\{0\}$, by repeating the previous reasoning with $\mO$ is the place of $\Om$ and $\Ga_* \cap \mO$ in the place of $\Ga_*$, we obtain once again the strict inequality $\E_\infty(u_*, \mO) < \E_\infty(u_* + \phi, \mO)$. The theorem ensues. \qed

\ms

The proof of Theorem \ref{existence_theorem} is more involved and requires some preparation. We begin with an elementary preliminary result.

\bl \label{lemma5} Let $\Om\sub \R^n$ be an open set and $\F \in C^2(\Om\by \R)$ a function satisfying \eqref{1.3}. Then
\begin{align}
\phantom{aaaaaaaa} & \sgn\big( \F(x,\xi) \big) \, =\, \sgn(\xi), \label{2.1}
\\
&  \frac{1}{c} |\xi| \, \geq\,  \big| \F(x,\xi)\big| \, \geq\, c|\xi|,  \label{2.2}
\\
& \big( |\F(x,\cdot)|^p \big)_{\xi\xi}(\xi) \,  \geq  \, 0, \ \ \ \ \ \text{ if }\ p\geq \frac{1}{c^3}+1,  \label{2.3}
\\
& \F_\xi (x,\xi)\, |\xi|\, \geq  \, c^2  |\F(x,\xi)|,  \label{2.4}
\end{align}
for all $(x,\xi) \in \Om \by \R$, where $c>0$ is the same constant as in \eqref{1.3}.
\el

\BPL \ref{lemma5}. We first note that \eqref{2.1} is obvious, while  \eqref{2.2} follows, by integration,
from $F(x, 0) = 0$ and the bounds $c\leq\F_\xi(x,\cdot)\leq1/c$ of \eqref{1.3}. For \eqref{2.3},
we differentiate in $\xi$ and use \eqref{1.3}. Then
\[
\begin{split}
 \big( |\F|^p \big)_{\xi\xi} \, &=\ p|\F|^{p-2}\Big(\F\,\F_{\xi\xi}\,+\, (p-1)(\F_\xi)^2 \Big) \\
& \geq\ p|\F|^{p-2}\Big( \!\!-\frac{1}{c}\,+\,(p-1)c^2\Big),
\end{split}
\]  
which establishes the desired bound. Finally, \eqref{2.4} is a consequence of \eqref{1.3}, which gives $|\xi|\F_\xi(x,\xi) \geq c|\xi|$, and of \eqref{2.2}, which gives $c|\xi|\geq c^2|\F(x,\xi)|$.    \qed
\ms

We will construct the solutions to our problem by approximation with minimisers of $L^p$ functionals. Therefore, we need to understand the behaviour of the latter.

\begin{proposition} \label{proposition6} Suppose that $\F \in C^2(\Om \times \R)$ satisfies \eqref{1.3}.
Then for any $p > c^{-3} +1$ there exists a minimiser $u_p$ of $\E_p$ over the space $ W^{2,p}_{u_0}(\Om)$.
Moreover, $u_p$ is a weak solution to the Dirichlet problem for the Euler-Lagrange equation associated with the functional $\E_p$:
\[
\left\{
\begin{array}{rl}
\De \Big( \big| \F(\cdot,\De u)\big|^{p-2}\, \F(\cdot,\De u)\,\F_\xi (\cdot,\De u) \Big) \,=\, 0,  &\ \ \text{ in }\Om,
\\
u \,=\, u_0,  &\ \ \text{ on }\p \Om,
\\ 
\D u \,=\, \D u_0, &\ \ \text{ on }\p \Om.
\end{array}
\right.
\]
Furthermore, there exist a (global) minimiser $u_\infty$ of the functional $\E_\infty(\cdot,\Om)$ over the space $\mathcal{W}^{2,\infty}_{u_0}(\Om)$ such that $\E_{p}(u_{p}) \larrow \E_\infty(u_\infty,\Om)$ as $p\ri \infty$.
Also, there exists a subsequence $(p_\ell)_1^\infty$ such that
\[
\left\{
\begin{array}{cl}
 u_{p_\ell} \larrow u_\infty, & \text{ in }C^{1}(\overline{\Om}),
\\
\D^2 u_{p_\ell} \weak \D^2 u_\infty, & \text{ in }L^q(\Om,\R^{n\by n}), \text{ for all }q \in (1,\infty),
\end{array}
\right.
\]
as $\ell \ri \infty$.
\end{proposition}

\BPP \ref{proposition6}. By \eqref{2.2}--\eqref{2.3} of Lemma \ref{lemma5}, for $p>c^{-3}+1$ the functional $\E_p$ is convex in $W^{2,p}_{u_0}(\Om)$ and 
\[
\E_p(u) \, \geq \, c(\mL^n(\Om))^{-1/p} \|\De u\|_{L^p(\Om)},
\]
for any $u\in W^{2,p}_{u_0}(\Om)$. Since $u-u_0 \in W^{2,p}_0(\Om)$, by the Calderon-Zygmund $L^p$ estimates (e.g.\ \cite{GT, GM}) and the Poincar\'e inequality we have
a positive constant $c_0 = c_0(p,\Om)$  that
\beq \label{L^p-estimate}
\| \De u \|_{L^p(\Om)}\, \geq\, c_0\, \| u \|_{W^{2,p}(\Om)} \, -\, (c_0 +1)\| u_0 \|_{W^{2,p}(\Om)}.
\eeq
Hence, $\E_p$ is coercive on $W^{2,p}_{u_0}(\Om)$ and by setting
\[
e_p\, :=\, \inf\big\{ \E_p(u) \, : \, u\in W^{2,p}_{u_0}(\Om)\big\}
\]
we also have the bound 
\[
0 \, \leq \, e_p\, \leq\, \E_p(u_0)\, < \, \infty
\]
because $u_0 \in W^{2,p}(\Om)$. By applying the direct method of the Calculus of Variations (e.g.\ \cite{D}), we deduce the existence of a global minimiser $u_p \in W^{2,p}_{u_0}(\Om)$. Further, $\E_p$ is Gateaux differentiable at the minimiser as a result of the bound 
\[
\big||\F(\cdot,\xi) |^{p-1}\F_\xi(\cdot,\xi) \big| \, \leq \, C |\xi|^{p-1}
\]
and well-known results (see e.g.\ \cite{D, GM}). 

Consider a family of minimisers $(u_p)_{p \ge p_0}$ where 
\[
p_0\, := \, \big\{ \text{integer part of } \max\{n, c^{-3}\} +1 \big\}
\]
and fix $k \in \N$. For any $p\geq k$, by \eqref{2.2}, H\"older's inequality and the minimality we have
\beq \label{2.5}
c\, \| \De u_p\|_{L^k(\Om)}\big(\mL^n(\Om)\big)^{-1/k} \,  \leq\, \E_k(u_p)\, \leq\, \E_p(u_p) \, \leq\, \E_p(u_0)\, \leq \, \E_\infty(u_0,\Om)
\eeq
and hence $(\De u_p)_{p \ge p_0}$ is bounded in $L^k(\Om)$. By the previous arguments and \eqref{L^p-estimate}, we conclude that $(u_p)_{p \ge p_0}$ is bounded in $W^{2,k}_{u_0}(\Om)$ for any $k\in \N$.
By a standard diagonal argument, weak compactness and the Morrey theorem, there exists \[
u_\infty \in \bigcap_{1<k<\infty}W_{u_0}^{2,\infty}(\Om)
\]
such that the desired convergences hold true along a subsequence as $p_\ell \ri \infty$.
When we pass to the limit as $\ell \ri \infty$ in \eqref{2.5}, the weak lower semicontinuity of the $L^k$ norm implies  
\[
\| \De u_\infty\|_{L^k(\Om)} \, \leq \,  \frac{(\mL^n(\Om))^{1/k}}{c} \,\E_\infty(u_0,\Om).
\] 
Letting $k \ri\infty$ we obtain $\De u_\infty \in L^\infty (\Om)$. Thus, $u_\infty \in \mathcal{W}^{2,\infty}_{u_0}(\Om)$, as desired. It remains to show the convergence of $\E_p(u_p)$ and minimality of $u_\infty$.

H\"older's inequality and minimality show that 
\[
\E_p(u_p) \,\le\, \E_p(u_q)\, \le \, \E_q(u_q), \ \text{ whenever } p \le q.
\]
Therefore, the limit $\underset{p \ri \infty}{\lim} \E_p(u_p)$ exists. Since $u_p - u_\infty \in W^{2,p}_0(\Om)$,
for any $\phi \in \mathcal{W}^{2,\infty}_0(\Om)$ the minimality and H\"older's inequality imply
\beq  \label{2.6}
\begin{split}
\E_\infty(u_\infty,\Om)\, &=\, \underset{k\ri \infty}{\lim} \, \E_k(u_\infty)
\\
 &\leq \, \underset{k\ri \infty}{\lim\inf}  \Big( \, \underset{\ell \ri \infty}{\lim\inf}\,  \E_k(u_{p_\ell}) \Big)
 \\
  &\leq \,  \lim_{p \ri \infty}  \E_p(u_p)
 \\
 &\leq \,  \underset{p\ri \infty}{\lim\sup}\,  \E_p(u_\infty + \phi) 
 \\
 & \leq \E_\infty(u_\infty +\phi ,\Om).
\end{split}
\eeq
Inequality \eqref{2.6} implies that $u_\infty$ is indeed a global minimiser of $\E_\infty (\cdot,\Om)$ over $ \mathcal{W}^{2,\infty}_{u_0}(\Om)$.
In addition, the choice $\phi = 0$ in \eqref{2.6} gives 
\[
\E_\infty(u_\infty, \Om) \, \leq \, \lim_{p\ri \infty} \E_p(u_p) \, \leq \, \E_\infty(u_\infty, \Om)
\]
and hence $\E_p(u_p) \larrow \E_\infty(u_\infty, \Om) $, as claimed.     \qed

\ms

The next result is an essential part of our constructions and this is the only point at which we make use of the $C^{1,1}$ boundary regularity of $\p\Om$ and the slightly higher regularity of the boundary condition $u_0 \in \mathcal{W}^{2,\infty}(\Om)$ (that is, that the Laplacian $\De u_0$ is continuous on $\overline{\Om}$). In fact, it suffices that $\De u_0$ be continuous near the boundary only. We are unsure whether it is possible to obtain the non-triviality of the limit without this condition. Certainly, the method we follow in Lemma \ref{lemma7} does not work.


Subsequently, we will be using the following symbolisation for the $r$-neighbourhood of the boundary $\p\Om$ in $\Om$:
\[
\Om_r \,:=\, \big\{ x\in \Om \, : \, \dist (x,\p \Om) <r \big\},
\]
for $r > 0$.

\bl[Improving the the boundary data] 
\label{lemma7}
Let $\Om \sub \R^n$ be a bounded open set with $C^{1,1}$ boundary and consider a function $u_0 \in \mathcal{W}^{2,\infty}(\Om)$ with $\De u_0 \in C(\overline{\Om})$. Then, for any $\epsilon > 0$ there exist a number $r=r(\e) > 0$ and
a function $w=w(\e) \in  \mathcal{W}_{u_0}^{2,\infty}(\Om)$ such that
 \[
 \| \De w \|_{L^\infty (\Om_r)}\, \leq\, \e.
 \]
\el

In other words, given any boundary condition $u_0  \in \mathcal{W}^{2,\infty}(\Om)$
such that $\De u_0$ is continuous up to the boundary, we can find another function $w$
in the same space with the same boundary data, the Laplacian of which is as small as desired in a
neighbourhood of the boundary. 

\BPL \ref{lemma7}. Since $\p\Om$ is $C^{1,1}$-regular, the result of the appendix establishes that the distance function $\dist(\cdot,\p\Om)$ belongs to $C^{1,1}(\overline{\Om_{2r_0}})$ for some $r_0>0$ small enough. 

Let $d $ be an extension of  $\dist(\cdot,\p\Om)$ from $\Om_{r_0}$ to $\overline{\Om}$ which is in the space $W^{2,\infty}(\Om)$. Extend $\De u_0$ continuously on $\R^n\set \Om$. Let $(\eta^\de)_{\de>0} \sub C^\infty_c(\R^n)$ be a standard mollifying family (as e.g.\ in \cite{E}). We set
\beq \label{2.7}
v_\de \, :=\, u_0\, -\, \frac{d^2}{2}\big(\eta^\de * \De u_0\big).
\eeq
Then $v_\de -u_0 \in \mathcal{W}^{2,\infty}_0(\Om)$, because $d=0$ on $\p\Om$ and also
\[
\begin{split}
\D v_\de \, &=\, \D u_0 \, - \, d\left\{\frac{d}{2} \D\big(\eta^\de * \De u_0\big) \,+\, \big(\eta^\de * \De u_0\big)\D d \right\},
\\
\D^2 v_\de \, &=\, \D^2 u_0 \, - \, \big(\D d \ot \D d\big)\big(\eta^\de * \De u_0\big)\, -\, d\bigg\{ \frac{d}{2} \D^2 \big(\eta^\de * \De u_0\big)
\\
& \ \  \ \ +\D\big(\eta^\de * \De u_0\big) \ot \D d \,+\, \D d \ot \D\big(\eta^\de * \De u_0\big) \,+\,
\big(\eta^\de * \De u_0\big)\D^2 d \bigg\} .
\end{split}
\]
By using that
\[
\tr(\D d \ot \D d) \, =\, |\D d|^2\, =\, 1 \ \ \text{ on $\Om_{r_0}$}, 
\]
for $0<r<r_0$ we deduce
\beq \label{2.8}
\begin{split}
\| \De v_\de \|_{L^\infty(\Om_r)}\, & \leq\, \left\| \De u_0 - \eta^\de *\De u_0 \right\|_{L^\infty(\Om_{r_0})} \, +\,   C_1 r \bigg( \big\| \De u_0 * \D^2 \eta^\de \big\|_{L^\infty(\R^n)} 
\\
& \ \ \ \ +\,  \big\| \De u_0 * \D \eta^\de \big\|_{L^\infty(\R^n)} \,+\, \big\| \De u_0  \big\|_{L^\infty(\R^n)}\bigg)
\end{split}
\eeq
for some constant $C_1 > 0$. Since 
\[
\eta^\de(x)\, =\, \de^{-n} \eta(|x|/ \de)
\]
for some fixed function $\eta \in C^\infty_c(\mB_1(0))$, by Young's inequality for convolutions we have the estimate
\beq  \label{2.9}
\ \ \  \big\| \De u_0 * \D^k \eta^\de \big\|_{L^\infty(\R^n)} \, \leq\, \frac{1}{\de^k} \, \big\| \De u_0  \big\|_{L^\infty(\R^n)} \, \|\D^k\eta\|_{L^1(\R^n)}, \ \ \ \ k=1,2.
\eeq
Hence, by invoking \eqref{2.9} we see that \eqref{2.8} gives 
\[
\begin{split}
\| \De v_\de \|_{L^\infty(\Om_r)}\, & \leq\, \left\| \De u_0 - \eta^\de *\De u_0 \right\|_{L^\infty(\Om_{r_0})}  +\,  C_2 r \bigg( \frac{1}{\de^2}\, +\, \frac{1}{\de}\, +\, 1  \bigg) \big\| \De u_0  \big\|_{L^\infty(\Om_{r_0})},
\end{split}
\]
where $C_2 = C_1 \|\eta\|_{W^{2,1}(\R^n)}$. By choosing 
\[
\de\, := \, r^{1/4}
\]
and also choosing $r_0$ sufficiently small, we obtain the desired statement as a consequence of the continuity of $\De u_0$.   \qed

\ms

Now we can show that the minimiser $u_\infty$ obtained in Proposition \ref{proposition6}  satisfies the desired formula of part \ref{(I)} in Theorem \ref{existence_theorem}. The rest of the proof is then not difficult.

\BPT \ref{existence_theorem}. Let us begin by setting
\[
e_\infty \, :=\,  \inf\big\{ \E_\infty(u, \Om) \, :\,  u \in \mathcal{W}^{2,\infty}_{u_0}(\Om)\big\}.
\]
If $e_\infty = 0$, then everything in Theorem \ref{existence_theorem} follows from Proposition \ref{proposition6}
or is trivial. (Note that in this case $e_p = 0$ for every $p$ and every minimiser of the corresponding
functionals is a harmonic function.) Therefore, we may assume that $e_\infty > 0$.

Since $e_p \larrow e_\infty$ as $p\ri \infty$ by Proposition \ref{proposition6}, it follows that $e_p>0$ for large $p$, say for $p \ge p_0$.
Then for all $p\geq p_0$, the formula in \ref{(IV)} gives rise to a measurable function $f_p : \Om \larrow \R$.
Then the Euler-Lagrange equation in Proposition \ref{proposition6} can be expressed in the form
\beq  \label{2.12}
\De f_p \, =\, 0, \ \ \text{ on }\Om.
\eeq
That is, $f_p$ is harmonic on $\Om$ and hence belongs to $C^\infty(\Om)$.

Let $p' =p/(p-1)$ be the conjugate exponent of $p \in (1,\infty)$. Then \eqref{1.3} implies
\[
\begin{split}
\left(\fint_\Om |f_p|^{p'}\right)^{1/p'} & =\, \frac{1}{e_p^{p - 1}} \bigg( \fint_\Om \Big| \F^{p-1} (\cdot,\De u_p)\F_\xi (\cdot,\De u_p) \Big|^{p'}\bigg)^{1/p'}
\\
& \leq \, \frac{1}{c \,e_p^{p - 1}} \left( \fint_\Om \Big| \F^{p-1} (\cdot,\De u_p)\Big|^{p/(p-1)}\right)^{(p-1)/p} \\
&= \, \frac{1}{c},
\end{split}
\]
which gives the following uniform $L^1$ bound of $(f_p)_{p \ge p_0} \sub C^\infty(\Om)$:
\beq  \label{2.13}
\|f_p\|_{L^1(\Om)} \, \leq \, \mL^n(\Om)\left( \fint_\Om |f_p|^{p'}\right)^{1/p'}\, \leq\, \frac{\mL^n(\Om)}{c}.
\eeq
By the mean value theorem for harmonic functions and by the standard interior derivative estimates (e.g.\ \cite{GT}) we have, for any $k\in \N\cup\{0\}$ and for any compactly contained $\mO \Subset \Om$, a constant $C = C(k,\mO,\Om)>0$ such that 
\[
\|\D^k f_p\|_{L^\infty(\mO)} \, \leq \, C \|f_p\|_{L^1(\Om)}.
\] 
Hence, the family $(f_p)_{p \ge p_0}$ is bounded (in the locally convex sense) in the topology of $C^\infty(\Om)$ and as a consequence there exist $f_\infty \in C^\infty(\Om)$ and a sequence $p_\ell \ri \infty$
such that $f_{p_\ell} \larrow f_\infty$ as $\ell \ri \infty$. From \eqref{2.12} it follows that  $f_\infty$ is harmonic:
\beq  \label{2.14}
\De f_\infty \, =\, 0, \ \ \text{ on }\Om.
\eeq
Fix $r>0$ and consider the inner $r$-neighbourhood $\Om_r$ of $\partial \Om$. Since $f_{p_\ell} \larrow f_\infty$ in $C\big(\overline{\Om \set \Om_r}\big)$, inequality \eqref{2.13} implies that
\[
\| f_\infty\|_{L^1(\Om\set \Om_r)} \, = \, \lim_{\ell \ri \infty} \, \|f_{p_\ell}\|_{L^1(\Om\set \Om_r)} \, \leq \,\frac{\mL^n(\Om)}{c}.
\]
Letting $r\ri 0$ we conclude that $f_\infty \in L^1(\Om)$ and
\beq \label{2.14a}
\| f_\infty\|_{L^1(\Om)} \,  \leq\, \frac{\mL^n(\Om)}{c}.
\eeq

We now show that $f_\infty \not\equiv 0$ using Lemma \ref{lemma7}. Fix $\epsilon>0$ small and let $r > 0$ and
$w \in \mathcal{W}^{2,\infty}_{u_0}(\Om)$ be as constructed in Lemma \ref{lemma7} with $| \De w| \leq \e$ on the inner zone $\Om_r \sub \Om$ of the boundary $\p\Om$. Since $u_p - w \in W^{2,p}_0(\Om)$, it is an admissible test function and by \eqref{2.12}, integration by parts gives
\[
\int_\Om f_p \, \De (u_p-w)\, =\, 0.
\]
Hence, by the above together with \eqref{def_f_p}, \eqref{1.3}, \eqref{2.1} and \eqref{2.4}, we obtain
\[
\begin{split}
\int_\Om f_p \, \De w \, & =\,  \int_\Om f_p \, \De u_p
\\
&=\, \frac{1}{e_p^{p - 1}} \int_\Om  \big| \F(\cdot,\De u_p)\big|^{p-2}\, \F(\cdot,\De u_p)\,\F_\xi (\cdot,\De u_p) \De u_p 
\\
&=\, \frac{1}{e_p^{p - 1}} \int_\Om  \big| \F(\cdot,\De u_p)\big|^{p-1}\, \F_\xi (\cdot,\De u_p) \, | \De u_p| 
\\
&\geq\,   \frac{c^2}{e_p^{p - 1}} \int_\Om  \big| \F(\cdot,\De u_p)\big|^{p}.
\end{split}
\]
Thus,
\beq  \label{2.15}
 \int_\Om f_p \, \De w \, \geq \, c^2\, \mL^n(\Om)\,  e_p.
\eeq
Now, we use \eqref{2.15}, \eqref{2.13} and Lemma \ref{lemma7} to estimate
\[
\begin{split}
c^2\, \mL^n(\Om)\,  e_p \, &\leq\, \int_{\Om_r}  f_p \, \De w \, +\, \int_{\Om\set \Om_r} f_p \, \De w 
\\
&\leq \, \e \|f_p\|_{L^1(\Om)} \, +\,  \int_{\Om\set \Om_r} f_p \, \De w 
\\
&\leq \,  \frac{\e \mL^n(\Om)}{c} \, +\,  \int_{\Om\set \Om_r} f_p \, \De w  .
\end{split}
\]
Recalling that $e_{p_\ell} \larrow e_\infty$ and also $f_{p_\ell} \larrow f_\infty$ in $C\big( \overline{\Om\set \Om_r} \big)$, we pass to the limit as $\ell \ri \infty$ to find that
\[
\int_{\Om\set \Om_r} f_\infty \, \De w \, \geq \, \mL^n(\Om) \Big(c^2\, e_\infty \, -\, \frac{\e}{c} \Big).
\]
By choosing $\e>0$ small enough, we deduce that
\[
\int_{\Om\set \Om_r} f_\infty \, \De w \, >\, 0,
\]
which implies $f_\infty \not\equiv 0$, as claimed.

Let us now define the map 
\[
\Phi \, : \,\ \Omega \times \R \larrow \Omega \times \R , \ \ \ \ \Phi(x, \xi) \,:=\, (x, \F(x, \xi)).
\]
Under the assumption \eqref{1.3}, this is a $C^2$ diffeomorphism.

Since the inverse function of $t\mapsto |t|^{p-2}t$ is given by $s\mapsto \sgn(s)|s|^{1/(p-1)}$, we may rewrite the formula \eqref{def_f_p} defining the harmonic function $f_p$, as
\[
\F(\cdot, \De u_p)\, =\, e_p |f_p|^{\frac{1}{p-1}} \big[ \F_\xi (\cdot,\De u_p)\big]^{-\frac{1}{p-1}} \sgn(f_p)
\]
or as
\beq \label{alternative_representation1}
\big(x, \De u_p(x)\big) = \Phi^{-1}\left( \cdot \,,\,  e_p |f_p|^{\frac{1}{p-1}} \big[ \F_\xi (\cdot,\De u_p)\big]^{-\frac{1}{p-1}} \sgn(f_p)\right)(x),
\eeq
for $x\in \Om \set\Ga_\infty$. On any compact set $K\sub \Om \set \Ga_\infty$, we have the uniform convergence
$f_{p_\ell} \larrow f_\infty$ as $\ell \ri \infty$, whereas $F_\xi$ is uniformly bounded from above and below by \eqref{1.3}. Hence, by restricting ourselves along the subsequence $p_\ell$ and letting $\ell \to \infty$ we obtain uniform convergence of the right-hand side of
\eqref{alternative_representation1} to $\Phi^{-1}\left(\cdot, e_\infty \sgn(f_\infty)\right)$ on $K$. But since we already know that $\De u_p \weak \De u_\infty$
weakly in $L^2(\Om)$, it follows that
\beq \label{alternative_representation2}
\big(x, \De u_\infty(x)\big) = \Phi^{-1}\Big(x, e_\infty \sgn\big(f_\infty(x)\big)\Big), \ \ \ x\in K.
\eeq
As a consequence,
\[
\F\big(x, \De u_\infty(x)\big) \, = \, e_\infty \sgn\big(f_\infty(x)\big), \ \ \ x\in K.
\]
Now let us recall that $\mL^n(\Gamma_\infty) = 0$. This is a consequence of general regularity results for nodal sets of solution to elliptic equations \cite{HS}. The statement of item \ref{(I)} then follows.

In order to prove item \ref{(II)}, we note that \eqref{alternative_representation2} implies that
$\De u_\infty \in C^2(\Om \set \Ga_\infty)$. The desired statement then follows from standard Schauder theory \cite{GT}.

For item \ref{(III)}, first recall the subsequential convergence of the $\E_p$-minimisers $(u_{p})_{1}^\infty$ of Proposition \ref{proposition6} along $(p_\ell)_1^\infty$ as $\ell\ri \infty$. We also have the desired respective convergence of the global infima $(e_{p})_{1}^\infty$ of the energies. 

The a.e.\ convergence of the Laplacians $(\De u_{p_\ell})_{\ell=1}^\infty$ follows from the fact that $\Om \set \Ga_\infty$ has full Lebesgue measure and that the sequence converges locally uniformly thereon. 

The strong convergence of the Laplacians $(\De u_{p_\ell})_{\ell=1}^\infty$ in $L^q(\Om)$ for all $q\in (1,\infty)$ is a consequence of the Vitaly convergence theorem (see e.g.\ \cite{FL}) and of the following facts:

\noi i) the weak convergence of the Laplacians over the same spaces, 

\noi ii) the a.e.\ convergence of the Laplacians on $\Om$, 

\noi iii) the boundedness of $\Om$, 

\noi iv) the $L^q$ equi-integrability estimate
\[
\| \De u_{p_\ell} \|_{L^q(E)}\, \leq\, \left(\sup_{\ell\in \N}\, \| \De u_{p_\ell} \|_{L^{q+1}(\Om)} \right) \big(\mL^n(E)\big)^{\frac{1}{q(q+1)}} 
\]
which holds true for any measurable subset $E\sub \Om$. 

Finally, the statement in part \ref{(IV)} has already been proven.    \qed

\ms

We conclude this section by noting that Corollary \ref{corollary4} is an immediate consequence
of Theorem \ref{uniqueness}, Theorem \ref{existence_theorem}, and the observation that
$\mL^n(\Gamma_\infty) = 0$ if $e_\infty >0$. On the other hand, $u_\infty$ is a harmonic function if $e_\infty=0$ and the result follows trivially.

\ms

\section*{Appendix: $C^{1,1}$ regularity of the distance function for $C^{1,1}$ domains} \label{section3}

Suppose that  $\Omega \sub \R^n$ is a bounded open set with  $C^{1, 1}$ boundary $\p\Om$. Let
\[
d \, \equiv \, \dist(\cdot,\p\Om)\ : \ \Omega \larrow \R
\]
symbolise the distance function to the boundary. For $r > 0$, let $\Omega_r$ denote again the inner $r$-neighbourhood of $\p\Om$ in $\Om$:
\[
\Om_r\,=\, \{x \in \Omega \,:\, d(x) < r\}. 
\]
In this appendix we establish that $d \in C^{1, 1}(\overline{\Omega_r})$ when $r$ is sufficiently small. This fact is
known, and the underlying tools are contained, e.g., in a paper by Feldman \cite{Fe}, but
we couldn't locate a precise reference in in the literature of the required form.
The $C^2$ regularity of the distance function for a $C^2$ boundary $\p\Om$ is a classical result, see e.g.\ \cite[Appendix 14.6]{GT}. On the other hand, the case of $C^1$ regularity of the distance function when the boundary $\p\Om$ is $C^1$ holds under the extra hypothesis that the distance is realised at \emph{one} point; see e.g.\ \cite{Fo}.

\ms

In order to prove the desired $C^{1,1}$ regularity of the distance function near $\p\Om$ when the boundary itself is a $C^{1,1}$ manifold (which we utilised in Lemma \ref{lemma7}), we first note the following fact: suppose that $r >0$ is such that $1/r$ is larger than the essential supremum of the curvature
of $\partial \Om$. If $x \in \Om$ and $y \in \partial \Om$ with $|x - y| = s \le r$ and such that
the tangent hyperplanes of $\partial \Om$ and $\partial \mB_s(x)$ coincide at $y$, then it follows that
$\mB_s(x) \sub \Om$ and $\partial \mB_s(x) \cap \partial \Om = \{y\}$. (This is easy to see when
$\partial \Om$ is $C^2$ regular and follows by an approximation of $\partial \Om$ with
$C^2$ manifolds otherwise.) Therefore, in the above situation, it follows that $d(x) = s$.
Moreover, for $x \in \Omega$ with $d(x) \le r$, it follows that there exists a \emph{unique}
point $y \in \partial \Om$  such that $|x - y| = d(x)$. Moreover, if $\nu$ denotes the outer
normal vector on $\partial \Omega$, then $x = y - d(x) \nu(y)$
and $\D  d(x) = - \nu(y)$.

\ms

Now fix $x_0 \in \partial \Omega$. Our aim is to prove $C^{1,1}$ regularity
of $d$ near $x_0$. To this end, we may assume without loss of generality that
there exist open sets $U \sub \R^n$ and $V \sub \R^{n - 1}$ such
that $x_0 \in U$ and $0 \in V$ and there exists a function $f \in C^{1, 1}(V)$
such that 
\[
\Omega \cap U \, =\, \Big\{(x', x_n) \in V \times \R\, :\,  x_n > f(x') \Big\} \cap U
\]
and such that $x_0 = (0, f(0))$ and $\D  f(0) = 0$. Define
\[
N(x')\, :=\, \frac{(-\D  f(x'), 1)}{\sqrt{1 + |\D  f(x')|^2}}, \quad x' \in V,
\]
so that $N(x') = -\nu(x', f(x'))$ for $x' \in V$. Note that this map is Lipschitz continuous.
We now define $\Psi : V \times \R \larrow \R^n$ by 
\[
\Psi(x', t)\, := \, (x', f(x')) \, +\, t\,N(x'). 
\]
Then $\Psi$ is $C^{0, 1}$ near $(0, 0)$.
Note also that $\Psi$ is injective in a sufficiently small neighbourhood of $(0, 0)$: if we had
$\Psi(x', s) = \Psi(y', t) =: z$ with $0 < s \le t \le r$, then it would follow that
$\partial \mB_t(z)$ and $\partial \Omega$ have the same tangent hyperplanes at $(y', f(y'))$.
By the above observations, this would imply that $\mB_t(z) \sub \Om$ and
$\partial \mB_t(z) \cap \partial \Om = \{(y',f(y'))\}$, and therefore $x' = y'$ and $s = t$.
So if $U$ and $V$ are chosen appropriately,
then $\Psi$ is a bijection between $V \times [0, r)$ and $U \cap \overline{\Omega}$. Moreover, we compute
\[
\D_j \Psi_i(x', t) \, =\,  \delta_{ij} \,+\,  t\, \D_j N_i (x'), \ \ \quad i, j = 1, \ldots, n - 1,
\]
and
\[
\D_j \Psi_n(x', t) \,=\, \D_j f(x') + t \,\D_j N_n(x'),\ \  \quad j = 1, \ldots, n - 1,
\]
while
\[
\D_t \Psi(x', t) \,=\, N(x').
\]
Since $\D  f(0) = 0$ and $N(0) = (0, \ldots, 0, 1)$, it follows that $\D \Psi$ is of full rank
in some neighbourhood of $(0, 0)$ and moreover, the inverse $(\D \Psi)^{-1}$ is
essentially bounded in this neighbourhood. That is, by making $V$ and $r$ smaller if necessary, without loss of generality we may assume that 
\[
\Psi^{-1} \,\in\, C^{0, 1}\big(U \cap \overline{\Omega};V \times [0, r)\big).
\]
Now note that
\[
\D  d(\Psi(x', t)) \,=\, N(x')
\]
whenever $t > 0$ is small enough. Hence if $\pi$ denotes the projection onto $\R^{n - 1} \times \{0\}$,
then we obtain the formula
\[
\D  d \, = \, N \circ \pi \circ \Psi^{-1}
\]
near $x_0$. The right-hand side is of class $C^{0, 1}$, and
thus $d$ is of class $C^{1, 1}$ near $x_0$. A compactness argument then
proves the above statement.

\ms
\ms

\noi \textbf{Acknowledgement.} N.K. would like to thank Craig Evans, Robert Jensen,  Jan Kristensen, Juan Manfredi, Giles Shaw and Tristan Pryer for inspiring scientific discussions on the topic of $L^\infty$ variational problems.

\ms

\noi \textbf{Funding:} N.K. has been partially financially supported by the EPSRC grant EP/N017412/1.

\noi \textbf{Conflict of interest:}
The authors declare that they have no conflict of interest.

\ms

\bibliographystyle{amsplain}

\end{document}